\documentclass[12pt]{article}

\usepackage{amsthm}

\theoremstyle{plain} 
\newtheorem{theorem}{Theorem}[section]
\newtheorem{thm} [theorem] {Théorème} 

\newtheorem{lm} [theorem]{Lemme} 

\newtheorem{prop}[theorem] {Proposition} 

\newtheorem{cor} [theorem] {Corollaire}

\newtheorem{rmq}[theorem] {Remarque}

\theoremstyle{definition} 

\usepackage[utf8]{inputenc}

\usepackage{amssymb,amsmath, amsfonts}

\usepackage{MnSymbol}

\usepackage{graphicx,color}

\author{ {St\'{e}phane DUGOWSON}\footnote{ISAE-Supmeca. Courriel : stephane.dugowson@isae-supmeca.fr}}
\date{22 février 2025}

\title{Un théorème brunnien pour les familles finies de variables aléatoires}

\begin{document}
\maketitle

\noindent\textbf{Abstract}. In 2014, during a study on the connectivity structures of quantum entanglement, I specifically introduced the notion of the connectivity structure of a family of random variables — a structure that expresses the dependency relations between the variables in question — and I stated the following proposition (\cite{Dugowson:20140718} p. 42), which can be described as Brunnian in reference to Hermann Brunn’s work on links \cite{Brunn:1892a} :
\emph{Every finite connectivity structure is that of a family of random variables}. At the time, however, I neglected to write down the proof of this assertion, merely providing an intuitive idea of it. The purpose of this article is to present such a proof.\\


\noindent\textbf{Résumé}. En 2014, à l'occasion d'un travail sur les structures connectives de l'intrication quantique \cite{Dugowson:20140718}, j'ai défini la notion de \emph{structure connective d'une famille de variables aléatoires} et j'ai énoncé la proposition suivante, qui peut être qualifiée de brunnienne en référence aux travaux de Hermann Brunn sur les entrelacs \cite{Brunn:1892a} : \emph{toute structure connective finie est celle d'une famille de variables aléatoires}. L'objet de cet article est d'en présenter une démonstration complète.\\

\noindent \textbf{Mots clés}  : Structures connectives, variables aléatoires, indépendance.

\section*{Introduction (English)} 

Without having formalized the notion of a connectivity structure, Hermann Brunn essentially stated in his 1892 paper that every finite (integral\footnote{Connectivity structures are said to be \emph{integral} when all singletons are connected.}) connectivity structure can be realized by a link embedded in the ambient space $\,\mathbf{R}^3$. While Brunn did indeed propose a construction of such a link based on what Rolfsen would call \emph{brunnian} links in 1976 \cite{Rolfsen:1976}, it was not until 1961, with the works of Debrunner \cite{Debrunner:19600416, Debrunner:1964}, that a rigorous proof was given for a similar statement (albeit in an ambient space of dimension $\geq 4$). Then, in 1985, Kanenobu \cite{Kanenobu:198504, Kanenobu:1986} provided a proof of the Brunnian claim in dimension $3$. As for the notion of connectivity structure and the category of connectivity spaces, they seem to have been explicitly introduced for the first time by Börger in 1981 \cite{Borger:1981, Borger:1983}, and we ourselves began a more in-depth study of them starting in the early 2000s \cite{Dugowson:200804, Dugowson:201012, Dugowson:201306, Dugowson:20180305}.

Whenever it is possible to associate a connectivity structure with the objects of a certain class, one can, in particular, raise the \emph{Brunnian} question of whether every connectivity structure of a certain type belongs to an object in that class. For example, in 2015 we established that any integral connectivity structure on a set $I$ (finite or otherwise) can be realized by a multiple relation \cite{Dugowson:20150505}. In the same vein, while working around the same time on quantum entanglement involving several particles, it was natural to wonder about the connectivity structure that expresses the dependence or independence of the random variables making up a family. Although this notion and the underlying idea for proving a Brunnian theorem for such random variable families were indeed laid out in \cite{Dugowson:20140718} (p.\ 42), I had not, at the time, written up a full proof. That proof, which turned out to be considerably more delicate than I had anticipated, is presented here in the case of finite integral connectivity structures.

\section*{Introduction (français)}

Sans avoir formalisé la notion de structure connective, Hermann Brunn avait en substance affirmé  dans son article publié en 1892 que toute structure connective intègre finie est celle d'un entrelacs plongé dans l'espace ambiant $\mathbf{R}^3$. Si Brunn avait effectivement proposé une construction d'un tel entrelacs fondée sur l'utilisation des entrelacs qui seront qualifiés de \emph{brunniens} par Rolfsen en 1976 \cite{Rolfsen:1976}, il faudra en fait attendre 1961 avec les travaux de Debrunner \cite{Debrunner:19600416, Debrunner:1964} pour avoir une démonstration rigoureuse d'un énoncé analogue, mais dans un espace ambiant de dimension $\geq 4$, puis 1985 avec ceux de Kanenobu \cite{Kanenobu:198504, Kanenobu:1986}  pour une démonstration de l'affirmation brunnienne en dimension $3$. La notion de structure connective et la catégorie des espaces connectifs semblent quant à elles avoir été explicitement considérées pour la première fois par Börger en 1981 \cite{Borger:1981, Borger:1983}, et nous en avons nous-même approfondi l'étude à partir du début des années 2000 \cite{Dugowson:200804, Dugowson:201012, Dugowson:201306, Dugowson:20180305}. Chaque fois qu'il est possible d'associer une structure connective  aux objets d'une certaine classe, on peut en particulier se poser la question \emph{brunnienne} de savoir si toute structure connective d'un certain type est celle d'un objet de cette classe. Nous avons par exemple établi en 2015 que toute structure connective intègre sur un ensemble $I$ (fini ou non) est celle d'une relation multiple \cite{Dugowson:20150505}. Dans le même ordre d'idées, travaillant à la même époque sur l'intrication quantique impliquant plusieurs particules, il était naturel de s'interroger sur la structure connective exprimant la dépendance ou l'indépendance de variables aléatoires constituant une famille. Si cette notion et l'idée devant présider à la démonstration d'un théorème brunnien pour de telles familles aléatoires ont bien été données dans \cite{Dugowson:20140718} (p. 42), je n'avais pas à l'époque rédigé une telle démonstration. Celle-ci, qui devait d'ailleurs s'avérer nettement plus délicate à établir que ce que je croyais, est donc proposée dans le présent article, dans le cas des structures connectives intègres finies.

\section{Notations et rappels}

Dans tout ce qui suit, sauf mention contraire, $n$ désigne un entier naturel non nul  et $I$ l'ensemble $\{1, 2, ..., n\}$. Pour tout ensemble $A$, nous désignons par  $\text{card}(A)$ ou $\vert A\vert$ son cardinal, et par $\mathcal{P}A$ l'ensemble de ses parties. Pour toute partie $J\in\mathcal{P}I$, nous notons $\neg J$ le complémentaire de $J$ dans $I$.

\paragraph{Rappels et notations relatives aux structures connectives.}

Pour l'ensemble des notions, défi\-nitions et notations relatives aux structures connectives, nous renvoyons le lecteur aux articles \cite{Dugowson:201012} et \cite{Dugowson:201306} (cette dernière référence comportant également quelques références historiques à ce sujet). Ainsi,  $\mathbf{Cnct}(I)$ désigne la catégorie (ordonnée par l'inclusion) des structures connectives \emph{intègres} sur $I$, c'est-à-dire des structures connectives pour lesquelles les singletons sont connexes.  Rappelons que toute structure connective $\mathcal{K}$ induit une structure connective $\mathcal{K}_J$ sur toute partie $J\in \mathcal{P}I$, à savoir  $\mathcal{K}_J=\mathcal{K}\cap \mathcal{P}J$. La structure connective intègre sur $I$ engendrée par un ensemble $\mathcal{J}$ de parties de $I$ est notée $[\mathcal{J}]$.
Pour tout $\mathcal{K}\in \mathbf{Cnct}(I)$, nous désignons par $\widetilde{\mathcal{K}}$ l'ensemble
$$\widetilde{\mathcal{K}}=
\{K\in \mathcal{K}, K\,\textrm{est\,irréductible, et \text{card}}(K)\geq 2\}.$$ 
L'ensemble $I$ étant fini, on a\footnote{Voir les propositions 5 et 16 de \cite{Dugowson:201012}.} $\mathcal{K}=[\widetilde{\mathcal{K}}]$ pour tout $\mathcal{K}\in\mathbf{Cnct}(I)$.   En outre, $\widetilde{\mathcal{K}_J}=\widetilde{\mathcal{K}}\cap \mathcal{P}J$.

\section{Dissociations ensemblistes, respect et contestations}

\subsection{Dissociations ensemblistes}

Nous appelons \emph{dissociation ensembliste}, ou simplement \emph{dissociation}, toute paire $\{J_1,J_2\}$ de parties de $I$ non vides et disjointes. Le \emph{domaine} $J$ d'une  dissociation est l'union, notée ${J=J_1\sqcup J_2}$, des deux parties qui la composent. Nous dirons également qu'une telle dissociation est \emph{une dissociation de $J$}. Les dissociations dont le domaine est $I$ seront dites \emph{globales}. Noter que le domaine d'une  dissociation a toujours au moins deux éléments. On note $\mathcal{D}$ l'ensemble des dissociations et, pour toute partie $J$ de $I$,  $\mathcal{D}(J)$ l'ensemble des \emph{dissociations de $J$ }  :
\begin{align*}
&\mathcal{D}=\{\{J_1,J_2\}\in\mathcal{P}\mathcal{P}I, (J_1\neq\emptyset\neq J_2)\wedge (J_1\cap J_2=\emptyset)\},\\
&\mathcal{D}(J) =\{\{J_1,J_2\}\in\mathcal{D}, J_1\sqcup J_2=J\}.
\end{align*}

Pour toute dissociation $\sigma$, on pourra également noter $\sigma_{(1)}$ et $\sigma_{(2)}$ les parties qui la constituent et  $\bar{\sigma}$ son domaine, de sorte que $\bar{\sigma}=\sigma_{(1)} \sqcup \sigma_{(2)}$. Soulignons que les parties $\sigma_{(1)}$ et $\sigma_{(2)}$ sont numérotées pour des raisons pratiques, mais qu'en fait il n'y a pas d'ordre entre elles (les échanger laisse la dissociation $\sigma$ inchangée).

\subsection{Contestations et respect des dissociations ensemblistes}\label{sec contestations dissociations ens}

Considérons une dissociation $\sigma\in\mathcal{D}$, dont nous noterons également $J=\bar{\sigma}$ le domaine et $J_i=\sigma_{(i)}$, où $i\in\{1,2\}$, les composantes. 
Nous dirons qu'une partie $A$ de $I$ \emph{conteste}, \emph{entrave} ou encore \emph{chevauche} $\sigma$, et nous écrirons $A\dashv \sigma$,  
si $A\subset J$ et $A\cap J_1\neq\emptyset\neq A\cap J_2$, autrement dit si $A\in\mathcal{P}J\setminus(\mathcal{P}J_1\cup\mathcal{P}J_2)$. \`{A} noter que, par définition, une dissociation est toujours contestée par son domaine\,: $\bar{\sigma}\dashv \sigma$.  \emph{A contrario}, si $A\in\mathcal{P}I$ ne conteste pas $\sigma$, nous écrirons $A\ndashv\sigma$ ou $A\rightspoon\sigma$, et nous dirons que $A$ \emph{valide} ou \emph{respecte} $\sigma$. 
Ainsi, pour toute partie $A$ de $I$, on a
$
(A\rightspoon\sigma)\Leftrightarrow (A\nsubset J) \vee (A\subset J_1) \vee (A\subset J_2)
$.
Pour tout ensemble de parties $\mathcal{A}\subset \mathcal{P}I$, posons
$$\overset{\mathcal{A}}{\dashv}\sigma=\{A\in\mathcal{A}, A\dashv\sigma\}=\{A\in\mathcal{A}\cap\mathcal{P}\bar{\sigma}, A\dashv\sigma\}.$$ 
Si $\overset{\mathcal{A}}{\dashv}\sigma\neq\emptyset$, nous dirons que $\mathcal{A}$ \emph{conteste} $\sigma$ et nous écrirons $\mathcal{A}\dashv\sigma$.  
Dans le cas contraire, nous dirons que $\mathcal{A}$ \emph{respecte} ou \emph{valide} $\sigma$, et nous écrirons $\mathcal{A}\ndashv \sigma$ ou $\mathcal{A}\rightspoon \sigma$, de sorte que l'on a les équivalences 
$$(\mathcal{A}\rightspoon \sigma)
\Leftrightarrow 
(\forall A\in\mathcal{A}, A\rightspoon\sigma)
\Leftrightarrow 
(\forall A\in\mathcal{A}\cap \mathcal{P}J, (A\subset J_1) \vee (A\subset J_2)).
$$

\subsection{Restrictions des dissociations}

La proposition suivante découle immédiatement des définitions. 

\begin{prop} \label{prop R - restriction d'une dissociation ensembliste}
Soient $L$ et $J$ deux parties de $I$ telles que $L\subset J$,  et soit $\sigma\in\mathcal{D}(J)$.  Si $L\dashv \sigma$ alors les traces sur $L$ de $\sigma_{(1)}$ et de $\sigma_{(2)}$ déterminent une dissociation $\sigma_{\vert L}$ de $L$, à savoir 
$
\sigma_{\vert L}=(L=(\sigma_{(1)}\cap L)\sqcup (\sigma_{(2)}\cap L))\in\mathcal{D}(L).
$
\end{prop}

La dissociation  $\sigma_{\vert L}$ ainsi induite sur $L$  sera appelée la \emph{trace} ou la \emph{restriction} à $L$ de la dissociation $\sigma\in\mathcal{D}(J)$.

\paragraph{Unions jointes.}
Les paires de parties non vides de $I$ qui \emph{ne sont pas} des dissociations joueront également un rôle important dans cet article. On appellera \emph{union jointe} l'union, notée $K_1 \cupdot K_2$, de deux telles parties $K_1$ et $K_2$\,:
 $$
(K= K_1 \cupdot K_2) \Leftrightarrow ((K = K_1 \cup K_2)\wedge (K_1\cap K_2\neq \emptyset)).
$$

\begin{lm} \label{lm dissociation d'une union jointe} Si $K=K_1\cupdot K_2$ et $\sigma\in \mathcal{D}(K)$, alors il existe 
$i\in\{1,2\}$ 
tel que $K_i\dashv \sigma$, 
et on a alors $\sigma_{\vert {K_i}}\in\mathcal{D}(K_i)$.
\end{lm}
\noindent \textit{Preuve}. Pour chaque $i\in\{1,2\}$ on a $K_i\subset K=\sigma_{(1)}\sqcup \sigma_{(2)}$, donc $K_i=( K_i\cap \sigma_{(1)}) \cup ( K_i\cap \sigma_{(2)})$, avec $( K_i\cap \sigma_{(1)}) \cap ( K_i\cap \sigma_{(2)})=\emptyset$. Il suffit donc de vérifier l'existence de $i\in\{1,2\}$ tel que $\forall j\in\{1,2\},  K_i\cap \sigma_j\neq \emptyset$. 
Si $K_2\cap\sigma_{(1)}=\emptyset$, on a $\sigma_{(1)}\subset (\sigma_{(1)}\cup \sigma_{(2)})\cap \neg K_2 = (K_1\cup K_2)\cap \neg K_2\subset K_1$, d'où $\emptyset\neq \sigma_{(1)}\cap K_1=\sigma_{(1)}$, et d'autre part $\sigma_{(2)}\cap K_2=(\sigma_{(1)}\cup \sigma_{(2)})\cap K_2= K_2$, donc $\sigma_{(2)}\supset K_2\supset K_1\cap K_2$ d'où $K_1\cap \sigma_{(2)}\supset K_2\cap K_1\cap \sigma_{(2)}=K_1\cap K_2\neq\emptyset$. De même, si $K_2\cap\sigma_{(2)}=\emptyset$, alors $\sigma_{(1)}\cap K_1 \neq\emptyset\neq \sigma_{(2)}\cap K_1$. Par conséquent, si la propriété $\forall j\in\{1,2\},  K_i\cap \sigma_j\neq \emptyset$ n'est pas satisfaite pour $i=2$, alors elle l'est nécessairement pour $i=1$, et c'est ce qu'il fallait vérifier.
\hfill $\square$

%
%

\subsection{Dissociations et connexité}   

\subsubsection{Dissociations globales adaptées aux composantes connexes}\label{subsub dissociations globales adapt composantes c}
Soit $\mathcal{K}$ une structure connective intègre sur $I$. Nous noterons $\mathcal{C}_\mathcal{K}$ l'ensemble des composantes connexes\footnote{Les composantes connexes d'un espace connectif sont les parties connexes maximales pour l'inclusion. Dans le cas d'un espace intègre, elles forment une partition de l'espace. Voir \cite{Dugowson:201012}.}   de l'espace connectif $(I,\mathcal{K})$ - souvent, pour simplifier l'expression, nous les désignerons aussi comme \emph{les composantes connexes de $\mathcal{K}$} - et nous noterons $\theta$ l'application $\mathcal{D}(\mathcal{C}_\mathcal{K})\rightarrow\mathcal{D}(I)$  qui à toute dissociation  $\gamma\in\mathcal{D}(\mathcal{C}_\mathcal{K})$ associe
la dissociation $\theta(\gamma)\in\mathcal{D}(I)$ 
définie pour chaque $i\in\{1,2\}$ par 
 $$\theta(\gamma)_{(i)}=\bigcup_{C\in \gamma_{(i)}} C.$$ Les dissociations de $I$ appartenant à $\theta(\mathcal{D}(\mathcal{C}_\mathcal{K}))$ seront dites \emph{adaptées aux composantes connexes de $\mathcal{K}$}.

\begin{prop}\label{prop dissociations globales} Les dissociations de $I$ adaptées aux composantes connexes de $\mathcal{K}$ sont celles que respectent tous les connexes $K\in\mathcal{K}$. Autrement dit, pour toute dissociation globale $\sigma\in\mathcal{D}(I)$, on a l'équivalence
$$ \mathcal{K}\rightspoon \sigma
\Leftrightarrow
\sigma\in \theta(\mathcal{D}(\mathcal{C}_\mathcal{K})).
$$
\end{prop}
\noindent \textit{Preuve}. Soit $\sigma\in\theta(\mathcal{D}(\mathcal{C}_\mathcal{K})$. Par construction, toute composante connexe $C\in \mathcal{C}_\mathcal{K}$ vérifie $C\subset\sigma_{(i)}$ pour un $i\in\{1,2\}$. Tout connexe $K\in\mathcal{K}$ étant inclus dans une telle composante $C$, on en déduit $K\rightspoon \sigma$, de sorte que $\mathcal{K}\rightspoon \sigma$. Réciproquement, si $\mathcal{K}$ respecte une dissociation globale $\sigma\in\mathcal{D}(I)$, toute composante connexe est en particulier incluse dans l'un des $\sigma_{(i)}$, de sorte que pour chaque $i\in\{1,2\}$ on a $\sigma_{(i)}=\bigcup_{C\in \gamma_{(i)}}C$, où $\gamma_{(i)}=\{C\in\mathcal{C}_\mathcal{K}, C\subset\sigma_{(i)}\}$, d'où $\sigma\in\theta(\mathcal{D}(\mathcal{C}_\mathcal{K}))$. \hfill $\square$

\subsubsection{Validation par les connexes irré\-ductibles}

\begin{prop}\label{prop validation par connexes irred}
Pour toute dissociation $\sigma=(J=J_1\sqcup J_2)\in \mathcal{D}$, on a les équivalences suivantes
$$
\mathcal{K}\rightspoon \sigma
\Leftrightarrow 
(\mathcal{K}_J=\mathcal{K}_{J_1}\cup\mathcal{K}_{J_2} )
\Leftrightarrow 
(\widetilde{\mathcal{K}_J}=\widetilde{\mathcal{K}_{J_1}}\cup\widetilde{\mathcal{K}_{J_2}} )
\Leftrightarrow 
\widetilde{\mathcal{K}}\rightspoon \sigma
$$
\end{prop}
 \noindent \textit{Preuve}. 
Si on a $\mathcal{K}\rightspoon (J=J_1\sqcup J_2)$, alors tout connexe $K\in\mathcal{K}\cap\mathcal{P}J$ est inclus dans l'un des $J_i$, autrement dit tout $K\in\mathcal{K}_J$ appartient à l'un des $\mathcal{K}_{J_i}$, et c'est en particulier le cas des connexes irréductibles ayant au moins deux éléments, de sorte que $\widetilde{\mathcal{K}}\rightspoon \sigma$. Réciproquement, si cette dernière propriété est satisfaite, alors $\mathcal{P}J_1 \cup \mathcal{P}J_2$ étant une structure connective sur $J$ qui contient $\widetilde{\mathcal{K}_J}$ on en déduit $[\widetilde{\mathcal{K}_J}]\subset \mathcal{P}J_1 \cup \mathcal{P}J_2$, d'où $\mathcal{K}_J\subset \mathcal{P}J_1 \cup \mathcal{P}J_2$ et donc $\mathcal{K}\rightspoon \sigma$. 
\hfill $\square$

\subsubsection{Caractérisation des connexes par dissociation}\label{sec connexes par dissociation}

Pour tout $\mathcal{A}\in\mathcal{P}\mathcal{P}I$, notons 
$\Gamma_\mathcal{A}
=\{B\in \mathcal{P}I, \forall \sigma\in\mathcal{D}(B), \mathcal{A} \dashv \sigma\}
$
 l'ensemble des parties $B$ de $I$ dont toute dissociation est contestée par au moins une partie appartenant à $\mathcal{A}$.

\begin{thm}\label{thm caractérisation connexe} Pour tout $\mathcal{K}\in \mathbf{Cnct}(I)$, on a $\Gamma_{\mathcal{K}}=\mathcal{K}$.
\end{thm} 
\noindent \textit{Preuve}. Si $C\in\mathcal{K}$ et $\sigma\in\mathcal{D}(C)$, alors $\mathcal{K}\ni C\dashv\sigma$, donc $C\in \Gamma_{\mathcal{K}}$. Réciproquement, soit  $A\in \Gamma_{\mathcal{K}}$. Comme tout espace connectif, l'espace induit sur $A$ se décompose en \emph{composante(s) connexe(s)} et si $A$ n'était pas  connexe, il aurait au moins deux composantes connexes. On pourrait dès lors former une dissociation  $\sigma=(A=A_1\sqcup A_2)$ en prenant pour chaque $A_i$ une union non vide de composantes connexes, de sorte que $\mathcal{K}_A=\mathcal{K}_{A_1}\cup \mathcal{K}_{A_2}$. Aucun connexe de $A$ ne contesterait alors $\sigma$, ce qui contredirait l'hypothèse de départ $A\in \Gamma_{\mathcal{K}}$.
\hfill $\square$

\begin{thm} Pour tout ensemble de parties $\mathcal{A}\subset \mathcal{P}I$, on a
$ \Gamma_\mathcal{A} = [\mathcal{A}]$.
\end{thm} 

\noindent \textit{Preuve}.
Vérifions d'abord que $\Gamma_\mathcal{A}$ est une structure connective intègre contenant $\mathcal{A}$ et que l'opérateur $\mathcal{A}\mapsto\Gamma_\mathcal{A}$ est croissant et qu'il laisse invariant les structures connectives intègres. L'égalité annoncée en découlera alors immédiatement, car $\mathcal{A}\subset [\mathcal{A}]$ $\Rightarrow$ $\Gamma_\mathcal{A}\subset\Gamma_{[\mathcal{A}]}=[\mathcal{A}]$ et réciproquement $\mathcal{A}\subset \Gamma_\mathcal{A}$ $ \Rightarrow$ $[\mathcal{A}]\subset [\Gamma_\mathcal{A}]=\Gamma_\mathcal{A}$.

\paragraph{$\Gamma_\mathcal{A}$ est une structure connective contenant $\mathcal{A}$.} 
Les singletons et la partie vide n'admettant pas de dissociation, ils sont nécessairement dans $\Gamma_\mathcal{A}$. De plus, étant données $B$ et $C$ appartenant à $\Gamma_\mathcal{A}$ et tels que $B\cap C\neq \emptyset$, toute dissociation $\sigma=J_1\sqcup J_2 \in\mathcal{D}(B\cupdot C)$ induit, d'après le lemme \ref{lm dissociation d'une union jointe}, une dissociation $\sigma_0$ de $B$ ou de $C$, nécessairement contestée par une partie $A\in \mathcal{A}$, qui à son tour conteste la dissociation $\sigma$ - par exemple, si la dissociation est induite sur $B$, on a  $\exists A\in\mathcal{A}, A \dashv \sigma_0=\sigma\cap B=(J_1\cap B, J_2\cap B)$  $\Rightarrow$ $A\cap J_1 \neq\emptyset\neq A\cap J_2$
- de sorte que $B\cup C \in\Gamma_\mathcal{A}$, ce qui prouve que $\Gamma_\mathcal{A}$ est bien une structure connective intègre. En outre, tout $A\in\mathcal{A}$ conteste toutes ses dissociations ensemblistes possibles, donc $\mathcal{A}\subset \Gamma_\mathcal{A}$.
\paragraph{$\mathcal{A}\mapsto\Gamma_\mathcal{A}$ est croissant.} Trivialement, car la communauté de ceux dont toute dissociation est contestée croît avec l'augmentation des instances  pouvant faire ces contestations. En effet, pour toute dissociation $\sigma\in\mathcal{D}$, on a : $$\mathcal{A}\subset\mathcal{B}\Rightarrow ((\exists A\in\mathcal{A}, A \dashv \sigma)\Rightarrow (\exists A\in\mathcal{B}, A \dashv \sigma)).$$
 
\paragraph{$\mathcal{A}\mapsto\Gamma_\mathcal{A}$ laisse les structures connectives intègres invariantes} C'est précisément ce qu'affirme le théorème \ref{thm caractérisation connexe}.
\hfill $\square$\\

Puisque $\Gamma_{\widetilde{\mathcal{K}}}=[\widetilde{\mathcal{K}}] =\mathcal{K}$, on en déduit immédiatement une caractérisation des connexes en terme de contestation de toute dissociation par les connexes irréductibles :

\begin{cor} \label{cor irred}
Pour tout structure connective intègre $\mathcal{K}$, on a 
$$
\mathcal{K}=
\{ A\in\mathcal{P}I, \forall\sigma\in\mathcal{D}(A), 
\exists K\in\widetilde{\mathcal{K}}, K \dashv \sigma
\}.
$$
\end{cor}

\section{Somme de structures connectives}\label{sec somme K}

On définit une application somme $\bigoplus:\mathbf{Cnct}(I)\times\mathbf{Cnct}(I) \rightarrow \mathbf{Cnct}(I)$  en posant pour tout couple $(\mathcal{K}_1,\mathcal{K}_2)$ de structures connectives intègres sur $I$ :
$
\mathcal{K}_1\oplus\mathcal{K}_2=[\mathcal{K}_1 \cup \mathcal{K}_2].
$ 
La catégorie $(\mathbf{Cnct}(I),\subset)$ étant un treillis complet, et la structure connective intègre $[\mathcal{A}]$ engendrée par un ensemble $\mathcal{A}$ de parties de $I$ étant la plus fine - \emph{i.e.} la plus petite pour l'ordre $\subset$ - des structures connectives intègres contenant $\mathcal{A}$ (voir \cite{Dugowson:201012}), l'opération $\oplus$ n'est rien d'autre que la borne supérieure : $\mathcal{K}_1\oplus\mathcal{K}_2=\mathcal{K}_1\vee\mathcal{K}_2$. La proposition suivante en découle immédiatement.

\begin{prop}\label{prop somme struct c asso commu} L'opération $\bigoplus:\mathbf{Cnct}(I)\times\mathbf{Cnct}(I) \rightarrow \mathbf{Cnct}(I)$ est \emph{associative}, \emph{commutative} et \emph{idempotente}. Elle admet la structure connective discrète intègre (pour laquelle seuls le vide et les singletons sont connexes) comme élément neutre, et la structure connective grossière (pour laquelle toutes les parties de $I$ sont connexes) comme élément absorbant.
\end{prop}
Nous étant placé dans le cadre d'un ensemble $I$ \emph{fini}, les structures connectives sont caractérisées par la donnée des connexes irréductibles, d'où la proposition suivante.
\begin{prop}\label{prop somme par les irred} Pour toute suite $\mathcal{K}_1, \mathcal{K}_2, ..., \mathcal{K}_r$ de structures connectives intègres, on a
$$
\bigoplus_{1\leq m\leq r}\mathcal{K}_m=
\bigvee_{1\leq m\leq r}{\mathcal{K}_m}=
[\bigcup_{1\leq m\leq r}{\mathcal{K}_m}]=
 [\bigcup_{1\leq m\leq r}\widetilde{\mathcal{K}_m}].
$$
\end{prop}
\noindent\textit{Preuve}. De $\widetilde{\mathcal{K}_m}\subset {\mathcal{K}_m}$, on tire $\bigcup_{1\leq m\leq r}\widetilde{\mathcal{K}_m}\subset \bigcup_{1\leq m\leq r}{\mathcal{K}_m}$ puis, par croissance de l'engendrement, $[\bigcup_{1\leq m\leq r}\widetilde{\mathcal{K}_m}]\subset [\bigcup_{1\leq m\leq r}{\mathcal{K}_m}]$. Réciproquement, $I$ étant fini, on a $\mathcal{K}_m=[\widetilde{\mathcal{K}_m}]\subset [\bigcup_{1\leq m\leq r}\widetilde{\mathcal{K}_m}]$, d'où $\bigcup_{1\leq m\leq r}{\mathcal{K}_m} \subset [\bigcup_{1\leq m\leq r}\widetilde{\mathcal{K}_m}]$, et donc, par croissance et idempotence de l'engendrement, $[\bigcup_{1\leq m\leq r}{\mathcal{K}_m}] \subset [\bigcup_{1\leq m\leq r}\widetilde{\mathcal{K}_m}]$. Mais puisque $[\bigcup_{1\leq m\leq r}{\mathcal{K}_m}]$ est la plus fine des structures connectives contenant tous les $\mathcal{K}_m$, elle peut aussi s'écrire $\bigvee_{1\leq m\leq r}{\mathcal{K}_m}$, c'est-à-dire $\bigoplus_{1\leq m\leq r}\mathcal{K}_m$.
\hfill $\square$

\section{Familles aléatoires}

\subsection{Définition des familles aléatoires sur $I$}

On appelle \emph{$I$-famille aléatoire} - ou \emph{famille aléatoire sur $I$}, ou simplement \emph{famille} - la donnée 
$$((\Omega, P), (X_1,..., X_n))$$ 
d'un espace de probabilité fini\footnote{Dont toute partie est supposée mesurable.} $(\Omega, P)$  et d'une famille de $n$ variables aléatoires (non nécessairement numériques\footnote{Bien que nous coderons souvent les valeurs de ces variables par des nombres, en particulier des nombres entiers écrits en base $2$, ces variables peuvent, comme dans la théorie de Shannon, être considérées comme essentiellement non numériques.}) $X_i:\Omega \rightarrow R_i\supset X_i(\Omega)$. L'ensemble $R_i$ sera appelé l'\emph{ensemble des valeurs aléatoires} ou  l'\emph{ensemble des valeurs possibles}  de $X_i$. 
Pour toute partie $J\subset I$, on pourra noter $X_J$ la variable aléatoire que constitue la famille $(X_j)_{j\in J}$, et $R_J$ l'ensemble produit $R_J = \prod_{j\in J} R_j \supset \prod_{j\in J} X_j(\Omega)\supset X_J(\Omega)$. 
On note $\mathcal{F}_I$, ou plus simplement $\mathcal{F}$, l'ensemble\footnote{En cas d'inquiétude quant au statut d'ensemble de $\mathcal{F}_I$ plutôt que de classe propre, on pourrait sans perte de généralité imposer que tous les ensembles $\Omega$ et tous les  $R_i=X_i(\Omega)$ soient inclus dans des ensembles de la forme $\{0,1\}^N$, avec $N$ des entiers qu'il serait même possible de borner en fonction de $n$.} de toutes les $I$-familles aléatoires. 
Pour toute famille aléatoire $\varphi\in \mathcal{F}_I$, on pourra noter $(\Omega^\varphi, P^\varphi)$ son espace de probabilité, $(X^\varphi_1,..., X^\varphi_n)$ ses variables aléatoires et, pour tout $i\in I$,  $R^\varphi_i$ les valeurs aléatoires de $X^\varphi_i$.

\subsection{Structure connective d'une $I$-famille aléatoire}

\subsubsection{Dissociations locales d'une famille aléatoire}\label{sec disso loc de phi}

Étant donnée une famille aléatoire $\varphi=((\Omega, P), (X_i:\Omega \rightarrow R_i)_{1\leq i \leq n})\in \mathcal{F}_I$ et $\sigma=(J=J_1 \sqcup J_2)\in\mathcal{D}(J)$ une dissociation de $\bar{\sigma}=J$ une partie de $I$, 
 nous dirons que  \emph{$\varphi$ respecte $\sigma$} ou que \emph{$\sigma$  dissocie $\varphi$ (localement)} ou encore que \emph{$\sigma$  dissocie $\varphi$ sur $J$} et nous écrirons
$\varphi\rightspoon\sigma$ ou  $\sigma\leftspoon \varphi$  si, pour toute famille $x_J \in R_J$,   $P(X_J=x_J)$ se factorise selon $J_1\sqcup J_2$ :
\[
(\sigma\leftspoon \varphi) 
\Leftrightarrow 
(\forall x_J\in R_J, P(X_J=x_J)=P(X_{J_1}=x_{J_1})P(X_{J_2}=x_{J_2})).
\]
Dans le cas particulier où $\bar{\sigma}=I$, autrement dit lorsque  $\sigma\in\mathcal{D}(I)$, une dissociation  $\sigma\leftspoon \varphi$ sera dite \emph{globale}.
Dans le cas où $\sigma\in\mathcal{D}(J)$ \emph{ne dissocie pas} $\varphi$ sur $\bar{\sigma}=J$, nous dirons que $\varphi$ \emph{conteste} $\sigma$, et nous écrirons $\varphi \dashv \sigma$, ou $\sigma \vdash \varphi$, ou $\sigma\nleftspoon \varphi$ ou encore $\varphi \nrightspoon \sigma$.

\begin{prop}\label{prop C localisation d'une dissociation familiale} Étant donné une partie $J\subset I$, une dissociation $\sigma\in\mathcal{D}(J)$ et une partie $A\subset J$ qui chevauche $\sigma$ (i.e. $A\dashv \sigma$), on a l'implication 
$$
(\varphi \rightspoon \sigma) \Rightarrow (\varphi \rightspoon \sigma_{\vert A}),
$$
ou, de façon équivalente, $(\varphi \dashv \sigma_{\vert A})\Rightarrow (\varphi \dashv \sigma).$
\end{prop}
 \noindent \textit{Preuve}. Posons $A_1=J_1\cap A$ et $A_2=J_2\cap A$, de sorte que $\sigma_{\vert A}=(A=A_1\sqcup A_2)$. On suppose que $\sigma$ dissocie sur $J$ la famille $\varphi$, et nous voulons vérifier que la dissociation a également lieu, \emph{a fortiori}, sur $A\subset J$, \emph{i.e.} $\sigma_{\vert A}\leftspoon \varphi$. Soit donc $x_A\in R_A$. Posons $M=J\setminus A$, $M_1=J_1\setminus A$ et $M_2=J_2\setminus A$. De
$$
(X_A=x_A)=\bigsqcup_{x_{M}\in R_{M}}(X_J=x_J),
$$
on tire, en sommant sur toutes les familles de valeurs $x_{M}\in R_{M}$,  
\begin{align*}
P(X_A=x_A)&=\sum_{x_{M}\in R_{M}}(X_J=x_J) 
=\sum_{x_{M_1}\in R_{M_1}}\sum_{x_{M_2}\in R_{M_2}}P(X_J=x_J) \\
&=\sum_{x_{M_1}\in R_{M_1}}\sum_{x_{M_2}\in R_{M_2}}P(X_{J_1}=x_{J_1})P(X_{J_2}=x_{J_2})\\
&=\sum_{x_{M_1}\in R_{M_1}}P(X_{J_1}=x_{J_1})\sum_{x_{M_2}\in R_{M_2}}P(X_{J_2}=x_{J_2})\\
&=P(X_{A_1}=x_{A_1})P(X_{A_2}=x_{A_2}),
\end{align*}
ce qu'il fallait vérifier.
\hfill $\square$

\subsubsection{Définition de la structure $\mathcal{K}_\varphi$}\label{sec def Kphi}

Conformément à la définition  donnée dans \cite{Dugowson:20140718}, p. 41 -  nous dirons que $J\subset I$ est \emph{non séparable pour $\varphi$} si celle-ci en conteste tout dissociation, autrement dit si $\forall \sigma\in\mathcal{D}(J), \varphi \dashv \sigma$. Notons, provisoirement, $\Gamma_\varphi$ l'ensemble des parties de $I$ non séparables pour $\varphi \in\mathcal{F}_I$ :
$$\Gamma_\varphi=\{J\in\mathcal{P}I,\forall\sigma\in\mathcal{D}(J), \varphi \dashv \sigma\}.$$
\begin{prop} L'ensemble  $\Gamma_\varphi$  est une structure connective intègre sur $I$. 
\end{prop}

\noindent\textit{Preuve}. Remarquons d'abord que les singletons sont trivialement non sépa\-rables pour $\varphi$. S'agissant de structures finies, il suffit alors de vérifier que l'union $K\cupdot L$ de \emph{deux} parties $K$ et $L$ d'intersection non vide, chacune étant non séparable pour $\varphi$,  est encore non séparable pour $\varphi$. On se donne donc deux telles parties.
Pour prouver que $K\cupdot L \in \Gamma_\varphi$, supposons par absurde qu'il existe $\sigma\in\mathcal{D}(K\cupdot L)$ telle que $\varphi\rightspoon \sigma$. D'après le lemme \ref{lm dissociation d'une union jointe}, cette dissociation $\sigma$ de $K\cupdot L$ conduirait par restriction à une dissociation ensembliste de $K$ ou de $L$, disons pour fixer les idées que $\sigma_{\vert K}\in\mathcal{D}(K)$. Mais d'après la proposition \ref{prop C localisation d'une dissociation familiale}, on devrait alors avoir $\varphi\rightspoon \sigma_{\vert K}$, ce qui contredirait l'hypothèse de non-séparabilité de $K$.
\hfill $\square$\\

Ayant ainsi vérifié que l'ensemble provisoirement noté $\Gamma_\varphi$ ci-dessus est bien une structure connective sur $I$, nous le noterons désormais $\mathcal{K}_\varphi$ et l'appel\-lerons \emph{la structure connective de $\varphi$}. Les $K\in \mathcal{K}_\varphi$ pourront être appelées les \emph{parties $\varphi$-connexes de $I$}.  L'application $\varphi\mapsto \mathcal{K}_\varphi$ qui à toute $I$-famille aléatoire associe sa structure connective sur $I$ pourra également être notée  $\kappa:\mathcal{F}_I  \mapsto\mathbf{Cnct}(I)$, de sorte que  $\kappa(\varphi)=\mathcal{K}_\varphi$. 

\subsubsection{Structure connective d'une famille restreinte}

Pour toute partie $J\subset I$, on définit de façon évidente l'application de restriction à $J$ des $I$-familles, $\mathcal{F}_I \ni \varphi \mapsto \varphi_{\vert J} \in \mathcal{F}_J$ en posant $\varphi_{\vert J}= ((\Omega,P), X_J)$ pour toute famille $\varphi=((\Omega,P), X_I=(X_1, ... X_n))$. 

\begin{prop}\label{prop Q struct conn de phi restreinte} Pour toute $\varphi\in \mathcal{F}_I$ et toute partie $J\subset I$, on a
$$ \kappa(\varphi_{\vert J})=
\mathcal{K}_{\varphi_{\vert J}}
=
{(\mathcal{K}_{\varphi})}_{\vert J}
=
{\mathcal{K}_\varphi}\cap\mathcal{P}J.
$$
\end{prop}
 \noindent \textit{Preuve}.  Pour toute partie $L\subset J$ et toute dissociation ensembliste $\sigma\in\mathcal{D}(L)$, la propriété $\varphi\dashv\sigma$ ne met en jeu que les variables $X_l$ avec $l\in L$, de sorte que l'on a trivialement
 $$(\varphi\dashv\sigma)\Leftrightarrow  (\varphi_{\vert J}\dashv\sigma) \Leftrightarrow  (\varphi_{\vert L}\dashv\sigma),$$ d'où l'on déduit que $L$ est $\varphi$-connexe si et seulement si elle est $\varphi_{\vert J}$-connexe, ce qui prouve la proposition énoncée. 
 \hfill $\square$
 
\subsubsection{Équivalence de la validation par $\varphi$ et par $\mathcal{K}_\varphi$}

Étant donnée $\varphi$ une famille aléatoire  sur $I$, notons $\mathcal{C}_\varphi = \mathcal{C}_{\mathcal{K}_\varphi}=\{C_1,...,C_p\}$ 
l'ensemble des composantes connexes de l'espace connectif $(I,\mathcal{K}_\varphi)$,
 numérotées de façon quelconque par les entiers de $1$ à $p=\mathrm{card}(\mathcal{C}_\varphi)$.
 Reprenant les notations de la section \ref{subsub dissociations globales adapt composantes c}, nous noterons $\theta(\mathcal{D}(\mathcal{C}_\varphi))$ l'ensemble des dissociations globales de $I$ adaptées à $\mathcal{C}_\varphi$, dont nous rappelons que, d'après la proposition \ref{prop disso glob adapt conn}, elles coïncident avec les $\sigma\in\mathcal{D}(I)$ telles que $\mathcal{K}_\varphi \rightspoon \sigma$.

 \begin{prop}[Indépendance des variables associées aux composantes connexes]\label{prop K ind V comp conn}   Les dissociations globales de $\varphi$, autrement dit les dissociations de $I$ qui sont respectées par $\varphi$, sont exactement les dissociations adaptées aux composantes connexes de $\mathcal{K}_\varphi$ :
$\{\sigma\in \mathcal{D}(I),  \varphi  \rightspoon \sigma\}=\theta(\mathcal{D}(\mathcal{C}_\varphi))$. 
Autrement dit, on a, pour toute dissociation $\sigma\in \mathcal{D}(I)$, l'équivalence
$$ 
(\varphi  \rightspoon \sigma) \Leftrightarrow (\mathcal{K}_\varphi \rightspoon \sigma).
$$
Il en découle en particulier que les variables aléatoires $X_{C_1}$, ..., $X_{C_p}$ associées aux composantes connexes sont mutuellement indépen\-dantes.
\end{prop}
\noindent \textit{Preuve}. Nous procédons en trois étapes, numérotées ci-après \textbf{1)} à \textbf{3)}. 

\noindent \textbf{1)} Commençons par vérifier que $\{\sigma\in\mathcal{D}(I), \varphi\rightspoon\sigma \}\subset \theta(\mathcal{D}(\mathcal{C}_\varphi))$, autrement dit que l'on a l'implication $\varphi  \rightspoon \sigma\Rightarrow\mathcal{K}_\varphi \rightspoon \sigma$. Soit donc $\sigma\in\mathcal{D}(I)$ telle que $\varphi  \rightspoon \sigma$. Si, par absurde, on avait $\mathcal{K}_\varphi \dashv \sigma$, alors il existerait un connexe $K\in\mathcal{K}_\varphi$ tel que $K\dashv\sigma$.  Mais, d'après la proposition \ref{prop C localisation d'une dissociation familiale}, on devrait alors avoir   $\varphi\rightspoon\sigma_{\vert K}\in\mathcal{D}(K)$, ce qui contredirait la $\varphi$-connexité de $K$. 

\noindent\textbf{2)} Montrons maintenant l'équivalence des deux propositions suivantes :
\begin{itemize}
\item a) toute dissociation issue des composantes connexes de $\mathcal{K}_\varphi$ est une 
dissociation globale de $\varphi$, autrement dit $\theta(\mathcal{D}(\mathcal{C}_\varphi))\subset \{\sigma\in\mathcal{D}(I), \varphi  \rightspoon \sigma\}$,
\item b) les variables aléatoires $X_{C_k}$, $1\leq k\leq p$, sont mutuellement indépen\-dantes.
\end{itemize}
Tout d'abord, il est clair que la proposition b) implique la proposition a), puisque l'indépendance mutuelle des $X_{C_k}$ permet d'écrire, pour toute dissociation $\sigma\in \theta(\mathcal{D}(\mathcal{C}_\varphi))$ et toute famille de valeurs $x_I\in R_I$,
\begin{align*}
P(X_I=x_I)&=\prod_{1\leq k \leq p} P(X_{C_k}=x_{C_k})\\
&=\prod_{k\in\sigma_{(1)}} P(X_{C_k}=x_{C_k})\prod_{k\in\sigma_{(2)}} P(X_{C_k}=x_{C_k})\\
&= P(X_{\sigma(1)}=x_{\sigma(1)})P(X_{\sigma(2)}=x_{\sigma(2)}).
\end{align*}
Réciproquement, supposons que la proposition a) soit valide. 
Pour tout ${k\in\{1,..., p-1\}}$, la dissociation  de $I$ définie par $\sigma^k=\theta(\gamma^k)$  où $\gamma^k$ est la dissociation de $\mathcal{C}_\varphi$ donnée par $\gamma^k=\{\{C_1,...,C_k\},\{C_{k+1}, ..., C_p\}\}$,  vérifie  ${\sigma^k\leftspoon \varphi}$. Puisque c'est en particulier le cas de $\sigma^1$, on en déduit d'abord que, pour toute famille de valeurs  $x_I\in R_I$, 
$$
P(X_I=x_I)=P(X_{C_1}=x_{C_1})P(\bigwedge_{k\geq 2} (X_{C_k}=x_{C_k})).
$$
Mais d'après la proposition \ref{prop C localisation d'une dissociation familiale}, la dissociation $\tau^2$ de $C_2\cup C_3\cup...\cup C_p$ 
obtenue par restriction de $\sigma^2$, 
\textit{i.e.}  $\tau^2={\sigma^2}_{\vert C_2\cup C_3\cup...\cup C_p}$, est également respectée par 
$\varphi$, de sorte que 
$$
P(\bigwedge_{k\geq 2} (X_{C_k}=x_{C_k}))=P(X_{C_2}=x_{C_2})P(\bigwedge_{k\geq 3} (X_{C_k}=x_{C_k}))
$$
et, en continuant ainsi de proche en proche jusque $k=p-1$, on obtient, du fait que $\varphi\rightspoon \tau^k={\sigma^k}_{\vert C_k\cup...\cup C_p}$, la factorisation complète qui assure l'indépendance mutuelle des $X_{C_k}$. \\

\noindent \textbf{3)} 
Prouvons enfin l'inclusion $\theta(\mathcal{D}(\mathcal{C}_\varphi))\subset\{\sigma\in\mathcal{D}(I), \varphi  \rightspoon \sigma\}$  par récurrence sur $p=\mathrm{card} (\mathcal{C}_\varphi)$. 
 Si $p=1$, on a $\mathcal{D}(\mathcal{C}_\varphi)=\emptyset$, 
 donc $\theta(\mathcal{D}(\mathcal{C}_\varphi))=\emptyset$, 
 et $I$ étant connexe dans ce cas, on a également $\{\sigma\in\mathcal{D}(I), \varphi  \rightspoon \sigma\} =\emptyset$. 
 Supposons maintenant que, quels que soient l'entier $n=\mathrm{card}(I)$ et la famille $\varphi\in\mathcal{F}_I$, 
 l'inclusion qui nous intéresse ait été établie pour toute valeur de $p\leq q$, où $q\in\{1,...,n-1\}$, 
 et montrons qu'elle est encore satisfaite pour $p=q+1$. 
 Puisque $\mathrm{card}(\mathcal{C}_\varphi)=q+1\geq 2$, $I$ n'est pas connexe, donc $\{\sigma\in\mathcal{D}(I), \varphi \rightspoon \sigma\}\neq\emptyset$. 
 Or, d'après le point \textbf{1)}, $\{\sigma\in\mathcal{D}(I), \varphi \rightspoon \sigma\}=\{\sigma\in\theta(\mathcal{D}(\mathcal{C}_\varphi)), \varphi \rightspoon \sigma\}$, 
 donc il existe au moins une dissociation $\sigma\in\theta(\mathcal{D}(\mathcal{C}_\varphi))$
  telle que $\varphi \rightspoon \sigma$. 
  On a alors, pour tout $x_I\in R_I$, 
  $$P(X_I=x_I)=P(X_{\sigma_{(1)}}=x_{\sigma_{(1)}})P(X_{\sigma_{(2)}}=x_{\sigma_{(2)}}),$$ mais puisque $\sigma$ est adaptée aux composantes connexes de $(I,\mathcal{K}_\varphi)$, l'ensemble des composantes connexes de $(\sigma_{(i)},\mathcal{K}_{\varphi_{\vert \sigma_{(i)}}})$ est, pour chaque $i\in\{1,2\}$, en conséquence de la proposition \ref{prop Q struct conn de phi restreinte} et de la définition des composantes connexes comme connexes maximaux,  $\mathcal{C}_{\varphi_{\vert \sigma_{(i)}}}=\mathcal{P}\sigma_{(i)} \cap \mathcal{C}_\varphi$, dont le cardinal est compris entre $1$ et $q$. L'hypothèse de récurrence s'applique alors à  la famille $\varphi_{\vert \sigma_{(i)}}$ et à toute dissociation de $\sigma_{(i)}$ adaptée aux composantes connexes de $\mathcal{K}_{\varphi_{\vert \sigma_{(i)}}}$, de sorte que 
 $$
 P(X_{\sigma_{(i)}}=x_{\sigma_{(i)}})=\prod_{C\in \mathcal{C}_\varphi, C\subset \sigma_{(i)}} P(X_C=x_C),
$$
d'où finalement $P(X_I=x_I)=\prod_{C\in \mathcal{C}_\varphi} P(X_C=x_C)$. Autrement dit, les variables aléatoires $X_C$ où $C$ décrit l'ensemble des composantes connexes de $\mathcal{K}_\varphi$ sont mutuellement indépendantes, et, d'après le point 2), on a bien $\theta(\mathcal{D}(\mathcal{C}_\varphi))\subset \{\sigma\in\mathcal{D}(I), \varphi  \rightspoon \sigma\}$.
 \hfill $\square$\\

\begin{cor}\label{cor K} Pour toute partie $A\subset I$, et pour toute dissociation $\sigma\in\mathcal{D}(A)$, on a
$ \varphi \rightspoon\sigma  \Leftrightarrow \mathcal{K}_\varphi \rightspoon\sigma$ (ou, de façon équivalente,
$\varphi \dashv \sigma  \Leftrightarrow \mathcal{K}_\varphi \dashv \sigma$.)
\end{cor}
 \noindent \textit{Preuve}.  On a en effet la suite d'équivalences suivantes 
$$
\mathcal{K}_\varphi \rightspoon\sigma
\Leftrightarrow
(\mathcal{K}_\varphi \cap \mathcal{P}A) \rightspoon\sigma
\Leftrightarrow
\mathcal{K}_{\varphi_{\vert A}}\rightspoon\sigma
\Leftrightarrow
\varphi_{\vert A}\rightspoon\sigma
\Leftrightarrow
\varphi\rightspoon\sigma,$$
où la première de ces équivalences découle de la définition de la validation d'une dissociation par un ensemble de parties (voir la section \ref{sec contestations dissociations ens}), la seconde de la proposition \ref{prop Q struct conn de phi restreinte}, la troisième de la proposition \ref{prop K ind V comp conn} appliquée à $\varphi_{\vert A}\in\mathcal{F}_A$ et à $\sigma$ considérée comme dissociation ``globale" de $A$, et la dernière de la définition d'une dissociation locale d'une famille aléatoire donnée en section \ref{sec disso loc de phi}.
  \hfill $\square$\\

\subsection{Produit tensoriel de deux $I$-familles aléatoires}

\subsubsection{Définition du produit tensoriel $\varphi\otimes\psi$}

On définit une application produit tensoriel $\bigotimes:\mathcal{F}_I\times\mathcal{F}_I \rightarrow \mathcal{F}_I$  en notant, pour tout couple $(\varphi,\psi)$ de familles aléatoires, $\varphi\otimes \psi$ la famille définie par
\begin{itemize}
\item $\Omega^{\varphi\otimes \psi}=\Omega^\varphi \times \Omega^\psi$, le produit cartésien des univers, 
\item $P^{\varphi\otimes \psi}=P^\varphi \otimes P^\psi$, le produit tensoriel des lois de probabilité, donné par $P^{\varphi\otimes \psi}(\omega, \xi)=P^\varphi(\omega)P^\psi(\xi)$ pour tous événements élémentaires $\omega\in \Omega^\varphi$ et $\xi\in \Omega^\psi$,
\item pour tout $i\in I$, $X^{\varphi\otimes \psi}_i:\Omega^{\varphi\otimes \psi} \rightarrow R^\varphi_i \times R^\psi_i$ est la variable aléatoire définie par $X^{\varphi\otimes \psi}_i=
X^\varphi_i \otimes X^\psi_i =
(\overline{X^{\varphi}_i},\overline{X^{\psi}_i})$, où  
 $\overline{X^{\varphi}_i} : \Omega^\varphi \times \Omega^\psi \rightarrow R^\varphi_i$
et 
$\overline{X^{\psi}_i} : \Omega^\varphi \times \Omega^\psi \rightarrow R^\psi_i$ 
sont définies pour tout $(\omega,\xi)\in \Omega^{\varphi\otimes \psi}$ par 
$$\overline{X^{\varphi}_i}(\omega,\xi) =X^\varphi_i(\omega)\,\,\mathrm{et}\,\,\overline{X^{\psi}_i}(\omega,\xi) =X^\psi_i(\xi),$$ autrement dit 
$X^{\varphi\otimes \psi}_i(\omega,\xi)=(X^{\varphi}_i(\omega),X^{\psi}_i(\xi))$.
\end{itemize}

\subsubsection{Indépendance des facteurs}

Considérons deux familles aléatoires notées $\varphi=((\Omega,P), (X_i:\Omega\rightarrow R_i)_{1\leq i \leq n})$ et $\psi=((\Lambda, Q), (Y_i:\Lambda\rightarrow S_i)_{1\leq i \leq n})$ , et soit $\varphi\otimes\psi=((\Theta,T),(Z_1,...,Z_n))$ leur produit tensoriel.

\begin{prop} \label{prop proba produit phi psi}
Quels que soient $A\subset \Omega$ et $B\subset\Lambda$, on a
$T(A\times B)=P(A)Q(B)$.
\end{prop}
 \noindent \textit{Preuve}. On a
$T(A\times B)=\sum_{\omega\in A,\, \xi\in B}T(\omega,\xi)
=\sum_{\omega\in A,\, \xi\in B}P(\omega)Q(\xi)$, d'où 
$T(A\times B)=\sum_{\omega\in A}P(\omega)\sum_{\xi\in B}Q(\xi)=P(A)Q(B)$. \hfill $\square$\\

\begin{prop}\label{prop independance XK et YL} Quelles que soient les parties $K\subset I$ et $L\subset I$, les variables $\overline{X}_K$ et $\overline{Y}_L$  sont indépendantes (sur leur univers commun $\Theta=\Omega\times \Lambda$).
\end{prop}
 \noindent \textit{Preuve}. Pour toutes familles  $x_K \in R_K$ et $y_L \in S_L$, on a
$ T((\overline{X}_K=x_K)\wedge(\overline{Y}_L=y_L)) 
=  T((X_K=x_K)\times(Y_L=y_L))$, expression qui,
d'après la proposition \ref{prop proba produit phi psi}, est égale à
$P(X_K=x_K)Q(Y_L=y_L)$, donc à $T(\overline{X}_K=x_K)T(\overline{Y}_L=y_L)$. 
\hfill $\square$\\

On en déduit le lemme suivant, valable pour tout couple $(\varphi,\psi)$ de telles familles aléatoires, et qui nous sera très utile dans la suite.

\begin{lm}  \label{lm dissociation entre phi et psi} Pour toute dissociation $\sigma\in\mathcal{D}$ on a
$$ {\varphi\otimes \psi}\rightspoon \sigma
\Leftrightarrow 
(\varphi\rightspoon \sigma) 
\wedge 
(\psi\rightspoon \sigma).$$
\end{lm}
\noindent\textit{Preuve}.  Posons $J=\bar{\sigma}$ et, pour chaque $i\in\{1,2\}$,  $J_i=\sigma_{(i)}$. 
Dire que ${\varphi\otimes \psi}\rightspoon \sigma$ équivaut à dire que pour toute famille de valeurs $(x,y)_J\in R_J \times S_J$   de la variable aléatoire $Z_J=(X\otimes Y)_J=(\overline{X},\overline{Y})_J$, on a
$$T\left(Z_J=(x,y)_J\right)=
T\left(Z_{J_1}=(x,y)_{J_1}\right)
\times
T\left(Z_{J_2}=(x,y)_{J_2}\right),$$
ce qui d'après la proposition \ref{prop proba produit phi psi} (ou la proposition \ref{prop independance XK et YL}) est équivalent à
$$
P(X_J=x_J)Q(Y_J=y_J)=P(X_{J_1}=x_{J_1})Q(Y_{J_1}=y_{J_1})P(X_{J_2}=x_{J_2})Q(Y_{J_2}=y_{J_2}).
$$
D'où, par les sommations d'usage, $P(X_J=x_J)=P(X_{J_1}=x_{J_1})P(X_{J_2}=x_{J_2})$ pour tous les $x_J$, d'où $\varphi\rightspoon \sigma$, et de même $Q(Y_J=y_J)=Q(Y_{J_1}=y_{J_1})Q(Y_{J_2}=y_{J_2})$ pour tous les $y_J$, d'où $\psi\rightspoon \sigma$. Mais puisque réciproquement ces relations entrainent la précédente, on en déduit l'équivalence annoncée.
\hfill $\square$

\subsubsection{Structure connective de $\varphi\otimes\psi$}\label{sec thm fond struct prod tens}

\begin{thm}\label{thm fondamental de la structure connective d'un produit tensoriel} Pour tout couple de familles aléatoires $(\varphi,\psi)\in (\mathcal{F}_I)^2$, on a $$\mathcal{K}_{\varphi\otimes\psi}=\mathcal{K}_{\varphi}\oplus\mathcal{K}_{\psi}.$$
\end{thm}
\noindent\textit{Preuve}. Pour toute partie $A\in\mathcal{P}I$, on a la suite d'équivalences suivantes :
\begin{flalign*}
A\in\mathcal{K}_{\varphi\otimes\psi} 
&\Leftrightarrow \forall \sigma\in\mathcal{D}(A), \varphi\otimes\psi\dashv\sigma
&&\text{(définition de }\mathcal{K}_\varphi,\text{ section \ref{sec def Kphi}})\\[1mm]
&\Leftrightarrow \forall \sigma\in\mathcal{D}(A), (\varphi\dashv\sigma)\vee (\psi\dashv\sigma)
&&\text{(lemme~\ref{lm dissociation entre phi et psi})}\\[1mm]
&\Leftrightarrow \forall \sigma\in\mathcal{D}(A), (\mathcal{K}_\varphi\dashv\sigma)
\vee(\mathcal{K}_\psi\dashv\sigma)
&&\text{(corollaire~\ref{cor K})}\\[1mm]
&\Leftrightarrow \forall \sigma\in\mathcal{D}(A), (\mathcal{K}_\varphi \cup \mathcal{K}_\psi) \dashv\sigma
&&\text{(par définition de }\mathcal{A}\dashv\sigma,\text{ section \ref{sec contestations dissociations ens})}\\[1mm]
&\Leftrightarrow A\in\Gamma_{(\mathcal{K}_\varphi \cup \mathcal{K}_\psi)}
&&\text{(par définition de }\Gamma_\mathcal{A}\text{, section \ref{sec connexes par dissociation})}\\[1mm]
&\Leftrightarrow A\in[{\mathcal{K}_\varphi \cup \mathcal{K}_\psi}]
&&\text{(théorème \ref{thm caractérisation connexe})}\\[1mm]
&\Leftrightarrow A\in {\mathcal{K}_\varphi \oplus\mathcal{K}_\psi}
&&\text{(définition de }{\mathcal{K}_1 \oplus\mathcal{K}_1,}\text{ section \ref{sec somme K}})\quad \hfill \square\\[1mm] 
\end{flalign*}

Pour tout ensemble de parties $\mathcal{A}\subset \mathcal{P}I$, la définition des parties irréductibles $K$ d'une structure connective $\mathcal{K}$ comme celles qui ne sont pas engendrées par les autres ($K\notin [\mathcal{K}\setminus \{K\}]$) entrainant le fait que $\widetilde{[\mathcal{A}]}\subset\mathcal{A}$, on déduit du théorème \ref{thm fondamental de la structure connective d'un produit tensoriel} le corollaire suivant.
\begin{cor} $\widetilde{\mathcal{K}_{\varphi\otimes\psi}}=\widetilde{\mathcal{K}_{\varphi}}\cup \widetilde{\mathcal{K}_{\psi}}$.
\end{cor}

Du corollaire précédent et des propositions \ref{prop validation par connexes irred} et \ref{prop somme par les irred}, on déduit alors  des équivalences figurant dans la démonstration  du théorème \ref{thm fondamental de la structure connective d'un produit tensoriel} la caractérisation suivante des $\varphi\otimes\psi$-connexes en termes de connexes irréductibles.
\begin{cor} Pour toute partie $A$ de $I$, on a l'équivalence
$$A\in\mathcal{K}_{\varphi\otimes\psi}
 \Leftrightarrow 
 \forall \sigma\in\mathcal{D}(A), (\widetilde{\mathcal{K}_\varphi}\dashv\sigma)
\vee(\widetilde{\mathcal{K}_\psi}\dashv\sigma).$$
\end{cor}

Du fait des propriétés de l'opération $\oplus$ sur les structures connectives (proposition \ref{prop somme struct c asso commu}), la proposition suivante découle également de façon immé\-diate du théorème \ref{thm fondamental de la structure connective d'un produit tensoriel}.
\begin{prop} Pour toutes $I$-familles aléatoires $\varphi$, $\psi$ et $\chi$, on a\\
$\mathcal{K}_{(\varphi\otimes\psi)\otimes \chi}=\mathcal{K}_{\varphi\otimes(\psi\otimes \chi)}$, 
$\mathcal{K}_{\varphi\otimes\psi}=\mathcal{K}_{\psi\otimes\varphi}$ et 
$\mathcal{K}_{\varphi\otimes\varphi}=\mathcal{K}_{\varphi}$.
\end{prop}

Il est clair qu'\emph{à des équi\-valences près} consistant en des isomorphismes entre espaces de probabilité, des changements de parenthésages et des permutations sur les familles de variables aléatoires et leurs valeurs, le produit tensoriel des $I$-familles aléatoires est \emph{associatif} et \emph{commutatif} : \emph{modulo} de telles équivalences, nous pourrions donc écrire $(\varphi\otimes\psi)\otimes \chi\simeq\varphi\otimes(\psi\otimes \chi)$ et $\varphi\otimes\psi\simeq \psi\otimes\varphi$,   d'où l'on aurait pu déduire 
 les strictes égalités correspondantes sur les structures connectives associées avant même d'avoir établi le théorème \ref{thm fondamental de la structure connective d'un produit tensoriel}. Ce dernier permet du moins de les établir rigoureusement et facilement sans avoir besoin de préciser formellement les équivalences en question. Soulignons que, par contre, $\varphi\otimes\varphi$ est en général une $I$-famille nettement plus complexe que $\varphi$.
Concluons cette section avec un autre corolaire immédiat du théorème \ref{thm fondamental de la structure connective d'un produit tensoriel} et de l'idempotence de $\oplus$.

\begin{cor}
Si  $\mathcal{K}_1$ et  $\mathcal{K}_2$ sont des structures connectives intègres telles que $\mathcal{K}_1\subset\mathcal{K}_\varphi\subset\mathcal{K}_2$ et $\mathcal{K}_1\subset\mathcal{K}_\psi\subset\mathcal{K}_2$, alors $\mathcal{K}_1\subset\mathcal{K}_{\varphi\otimes\psi}\subset\mathcal{K}_2$.
\end{cor}


\section{Familles brunniennes}

\subsection{Famille discrète canonique}

On appellera \emph{famille discrète canonique sur $I$} la $I$-famille aléatoire $\beta_d$ définie par $\beta_d=((\Omega_d, P_d),(B^d_1,...,B^d_n))$ avec
\begin{itemize}
\item $\Omega_d=\{0\}$, de sorte que $P_d(0)=1$,
\item $\forall i\in I, B^d_i=\text{Id}_{\Omega_d}$.
\end{itemize}
Bien que toutes les ``variables" de cette famille soient la même, elles sont mutuellement indépendantes car elles sont constantes (une variable constante est indépendante d'elle-même). Il en découle que $\mathcal{K}_{\beta_d}=\mathcal{K}_d$, la structure connective intègre discrète sur $I$ (\emph{i.e.} celle pour laquelle les seuls connexes sont $\emptyset$ et les singletons $\{i\}$, avec $i\in I$). 
Remarquons au passage qu'il existe une infinité d'autres familles aléatoires constituées de variables mutuellement indépendantes, en particulier toutes celles de la forme 
$$\varphi=((\Omega_\varphi=\prod_{i\in I} \Omega_i, \bigotimes_{i\in I} P_i),(X_i\circ \pi_i)_{i\in I})$$ où chaque $X_i$ est une variable aléatoire sur un espace de probabilité fini $(\Omega_i, P_i)$ et où $\pi_i:\Omega_\varphi \rightarrow \Omega_i$ désigne la $i$-ème projection. 

\subsection{Famille $M$-brunnienne canonique}

Soit $M=\{i_1,..., i_m\}\subset I$, avec $m=\text{card}(M)\geq 2$, les $i_k$ étant rangés par ordre croissant. On appellera \emph{famille $M$-brunnienne canonique sur $I$} la $I$-famille aléatoire $\beta_M=((\Omega_M, P_M),(B^M_1,...,B^M_n))$ définie par $\Omega_M=\{0,1\}^{m-1}$ et, pour tout $\omega=(\omega_1,...,\omega_{m-1})\in \Omega_M$,

\begin{itemize}
\item $P_M(\omega)=\dfrac{1}{2^{m-1}}$ (autrement dit $P_M$ est la probabilité uniforme sur $\Omega_M$),
\item $\forall k\in\{i_1,...,i_{m-1}\}$, $B^M_k(\omega)=\pi_k(\omega)=\omega_k\in\{0,1\}$,
\item $B^M_{i_m}(\omega)=(\sum_{k\in\{i_1,...,i_{m-1}\}} \omega_k )(\text{mod}\,2)\in\{0,1\}$, 
\item $\forall k\in I\setminus M, B^M_k(\omega)=0$.
\end{itemize}

\begin{prop}\label{prop famille brunnienne} Outre le vide et les singletons, la structure connective de $\beta_M$ admet $M\subset I$ pour seul connexe. Ainsi, $\mathcal{K}_{\beta_M}=\{\emptyset, \{1\}, ..., \{n\}, M\}$ ou, de façon équivalente $\widetilde{\mathcal{K}_{\beta_M}}=\{M\}$, ou encore :
$$
\kappa({\beta_M})=[\{M\}].
$$
\end{prop}
\noindent \textit{Preuve}. 
Si $J\subset I$ a au moins deux éléments dont au moins un $k\in I\setminus M$, alors $B^M_k$, étant constante, se dissocie de la famille $B^M_{J\setminus \{k\}}$, de sorte que $J\notin\mathcal{K}_{\beta_M}$. Par conséquent, toute partie $K\in \mathcal{K}_{\beta_M}$ qui a au moins deux éléments est nécessairement incluse dans $M$. Il est clair que $M$ lui-même est connexe : du fait de la relation entre les $m$ variables aléatoires $B^M_k$ pour $k\in\{i_1,..., i_m\}$, qui détermine la valeur de l'une quelconque d'entre elles en fonction des valeurs prises par les $m-1$ autres, la famille $B^M_M$ ne saurait être dissociée. Soit alors $J\subsetneqq M$ telle que $\text{card}(J)\geq 2$ (ce qui suppose $n\geq m\geq 3$). Si $i_m\notin J$, les variables $(B^M_j)_{j\in J}$ sont mutuellement indépendantes de façon évidente, par leur construction même, donc on a encore $J\notin \mathcal{K}_{\beta_M}$. Supposons donc $i_m\in J$, de sorte que $\emptyset \subsetneqq L \subsetneqq J \subsetneqq M \subset I$, où $L=J\setminus \{i_m\}$.  Quitte à ré-indexer les éléments de $I$ pour simplifier l'expression des éléments à considérer, on peut supposer sans restreindre la généralité du raisonnement qui va suivre que l'on a $i_k=k$ lorsque $1\leq k\leq m$, autrement dit $M=\{1, 2, ..., m\}$, et que $J=L\sqcup \{m\}$ avec $L=\{1,..., l\}$, le cardinal $l$ de $L$ vérifiant donc $1\leq l = \text{card}(J)-1\leq m -2$.  Soit alors $(b_J)\simeq (b_L,b_m)\in R_L \times \{0,1\} \simeq R_J = \{0,1\}^J$ une $J$-famille\footnote{Les équivalences $\simeq$ ne sont ici que de simples re-parenthésages.} de valeurs pour $B^M_J$. L'événement $(B^M_J=b_J)$ est l'ensemble des événements élémentaires $\omega=(\omega_1,...,\omega_{m-1})$ vérifiant
\begin{itemize}
\item $\forall k\in L, \omega_k=b_k$,
\item $\forall k\in \{l+1,..., m-2\}$, $\omega_k\in\{0,1\}$, soit $(m-2-l)$ degrés de liberté (éventuellement aucun degré de liberté si $l=m-2$, auquel cas $\{l+1,..., m-2\}=\emptyset$),
\item $\omega_{m-1}=b_m+(b_1+...+b_l)+(\omega_{l+1}+...\omega_{m-2}) (\text{mod}\,2)$,
\end{itemize}
de sorte que $\vert (B^M_J=b_J)\vert=2^{m-2-l}$, d'où (la probabilité $P_M$ étant uniforme  sur $\Omega_M$), $P_M(B^M_J=b_J)=2^{m-2-l} \times \dfrac{1}{2^{m-1}} = \dfrac{1}{2^{l+1}}$. Or, $P_M(B^M_1=b_1)=...=P_M(B^M_l=b_l)=P_M(B^M_m=b_m)=\dfrac{1}{2}$ d'où $P_M(B^M_J=b_J)=\prod_{j\in J}P_M(B^M_j=b_j)$, et donc, ici encore, $J\notin\mathcal{K}_{\beta_M}$. Finalement, comme annoncé, la seule partie de $I$ ayant au moins deux éléments et qui soit dans $\mathcal{K}_{\beta_M}$ est $M$.
 \hfill $\square$

\section{Construction d'une famille aléatoire de structure connective donnée}

\begin{thm}\label{thm brunnien} Pour toute structure connective intègre $\mathcal{K}\in \mathbf{Cnct}(I)$, il existe une famille aléatoire $\varphi\in\mathcal{F}_I$ telle que $\mathcal{K}_\varphi=\mathcal{K}$.
\end{thm}
\noindent \textit{Preuve}. Si $\mathcal{K}$ est la structure discrète, la famille discrète $\beta_d$ répond à la question. Sinon, l'ensemble $\widetilde{\mathcal{K}}=\{K_1,..., K_r\}$ est non vide, et le produit tensoriel des familles brunniennes canoniques associées aux $K\in \widetilde{\mathcal{K}}$ vérifie la propriété souhaitée en vertu du théorème \ref{thm fondamental de la structure connective d'un produit tensoriel}, de la proposition \ref{prop famille brunnienne} et de la proposition \ref{prop somme par les irred}  :
$$
\kappa(\bigotimes_{K\in \widetilde{\mathcal{K}}}\beta_K)
=
\bigoplus_{K\in \widetilde{\mathcal{K}}} \kappa(\beta_K)
=
\bigoplus_{K\in \widetilde{\mathcal{K}}} [\{K\}]
=
[\widetilde{\mathcal{K}}]=\mathcal{K}.
$$
\hfill $\square$

\begin{rmq} Bien que la construction ci-dessus dépende du choix de l'ordre dans lequel les connexes irréductibles $K_1, K_2,..., K_r$ sont rangés et de l'ordre dans lequel le produit tensoriel des familles brunniennes associées est effectué, les familles finalement obtenues devraient être équivalentes entre elles en un sens assez naturel que nous laissons au lecteur intéressé le soin de préciser (voir aussi la remarque similaire faite en section \ref{sec thm fond struct prod tens}).
\end{rmq}

\begin{rmq}\label{rmq trouver varphi optimal}
La preuve du théorème \ref{thm brunnien} donne un \emph{procédé universel explicite  de construction} d'une famille $\varphi$ telle que $\mathcal{K}_\varphi=\kappa(\varphi)=\mathcal{K}$. Soulignons toutefois que la famille ainsi obtenue n'est pas nécessairement optimale, au sens où d'autres familles sur un univers $\Omega$ de cardinal inférieur pourraient parfois être trouvés. Un exemple flagrant de cela est donné par la structure connective grossière sur $I$, à savoir $\mathcal{K}^g=\mathcal{P}I$, pour laquelle toute partie de $I$ est connexe. 
On a $\widetilde{\mathcal{K}^g}=\{\{i,j\}\in\mathcal{P}I, i\neq j\}$, d'où $\vert \widetilde{\mathcal{K}^g}\vert = \dfrac{n(n-1)}{2}$. 
Comme pour chaque $K\in\widetilde{\mathcal{K}^g}$ on a $\vert \Omega_K \vert = 2$, la construction proposée dans la démonstration du théorème \ref{thm brunnien} conduit à une famille $\varphi=\bigotimes_{K\in\widetilde{\mathcal{K}^g}} \beta_K$ dont l'univers possède $2^{\frac{n(n-1)}{2}}$ éléments. 
Des familles bien plus simples pour la même structure connective sont celles de la forme  $\varphi=((\Omega, P), (X, X,..., X))$ où toutes les variables sont non constantes et égales entre elles, avec par exemple $\Omega=\{0,1\}$ muni de la probabilité uniforme et $X=\text{Id}_{\Omega}$.
\end{rmq}

\begin{cor} Pour toute structure connective intègre $\mathcal{K}$, il existe une infinité de familles aléatoires $\varphi$ telles que $\kappa(\varphi)=\mathcal{K}$.
\end{cor}
\noindent \textit{Preuve}. On a déjà vu que c'était le cas pour la structure connective intègre discrète. Pour toute autre structure $\mathcal{K}$, la famille $\varphi$ construite dans la preuve du théorème \ref{thm brunnien} a un univers de cardinal au moins $2$, de sorte que les familles de la suite $\varphi$, $\varphi\otimes \varphi$, $\varphi\otimes\varphi\otimes \varphi$, etc... sont deux à deux distinctes : elles fournissent un ensemble infini de familles de même structure connective $\mathcal{K}$.\hfill $\square$

\begin{cor} Il existe une application $\wedge:\mathcal{F}_I\times \mathcal{F}_I \rightarrow \mathcal{F}_I $ telle que pour tout couple $(\varphi,\psi)\in {\mathcal{F}_I}^2$ on ait $\mathcal{K}_{\varphi\wedge\psi}=\mathcal{K}_\varphi \cap \mathcal{K}_\psi$.
\end{cor}
\noindent \textit{Preuve}. En choisissant un procédé universel pour ordonner les connexes irréductibles de toute structure connective sur $I$ - par exemple en utilisant l'ordre lexicographique induit par celui de $I$ - la construction du théorème \ref{thm brunnien} appliquée à la structure connective  $\mathcal{K}_\varphi \cap \mathcal{K}_\psi$ fournit un procédé universel de construction de $\varphi\wedge\psi$. Cela dit, il pourrait être intéressant de chercher un $\varphi\wedge\psi$ qui soit optimal (voir la remarque \ref{rmq trouver varphi optimal}).\hfill $\square$

\section{Conclusion} Cette étude est limitée au cas des familles finies de variables aléatoires, il pourrait être intéressant de chercher à l'étendre au cas infini. Par ailleurs, la question reste ouverte de la recherche de familles optimales (\textit{i.e.} sur un univers minimal)  de structure connective donnée. Il serait également intéressant d'établir de nouvelles propriétés brunniennes dans d'autres domaines, et de chercher à préciser les relations pouvant exister entre elles.


\bibliographystyle{plain}


\begin{thebibliography}{10}

\bibitem{Borger:1981}
Reinhard B\"orger.
\newblock {\em Kategorielle Beschreibungen von Zusammenhangsbegriffen}.
\newblock PhD thesis, Fernuniversit\"at, Hagen, 1981.

\bibitem{Borger:1983}
Reinhard B\"orger.
\newblock Connectivity spaces and component categories.
\newblock In {\em Categorical topology, International Conference on Categorical
  Topology (1983)}, Berlin, 1984. Heldermann.

\bibitem{Brunn:1892a}
Hermann Brunn.
\newblock Ueber verkettung.
\newblock {\em Sitzungsberichte der Bayerische Akad. Wiss., MathPhys. Klasse},
  22:77--99, 1892.

\bibitem{Debrunner:19600416}
Hans Debrunner.
\newblock {Links of Brunnian type}.
\newblock {\em Duke Math. J.}, 28:17--23, 1961.

\bibitem{Debrunner:1964}
Hans Debrunner.
\newblock {\"Uber den Zerfall von Verkettungen.}
\newblock {\em Mathematische Zeitschrift}, 85:154--168, 1964.
\newblock {http://www.digizeitschriften.de}.

\bibitem{Dugowson:20180305}
St\'{e}phane Dugowson.
\newblock {Toposes of connectivity spaces. Morita equivalences with topological
  spaces and partially ordered sets in the finite case}, 5 mars 2018.
\newblock \\\texttt{https://hal.archives-ouvertes.fr/hal-01722695}.

\bibitem{Dugowson:200804}
Stéphane Dugowson.
\newblock The connectivity order of links.
\newblock Article accessible en ligne sur
  {http://hal.archives-ouvertes.fr/hal-00275717/en/}, 2008.

\bibitem{Dugowson:201012}
Stéphane Dugowson.
\newblock On connectivity spaces.
\newblock {\em Cahiers de {T}opologie et {G}éométrie {D}ifférentielle
  {C}atégoriques}, LI(4):282--315, 2010.
\newblock http://hal.archives-ouvertes.fr/hal-00446998/fr.

\bibitem{Dugowson:201306}
Stéphane Dugowson.
\newblock Espaces connectifs : représentations, feuilletages, ordres,
  difféologies.
\newblock {\em {Cahiers de Topologie et Géométrie Différentielle
  Catégoriques}}, LIV, 2013.
\newblock \\\texttt{https://hal.archives-ouvertes.fr/hal-01386249}.

\bibitem{Dugowson:20140718}
Stéphane Dugowson.
\newblock {Structures connectives de l'intrication quantique}, 18 juillet 2014.
\newblock \texttt{https://hal.archives-ouvertes.fr/hal-01025949}.

\bibitem{Dugowson:20150505}
Stéphane Dugowson.
\newblock Structure connective des relations multiples, 5 mai 2015.
\newblock \texttt{https://hal.archives-ouvertes.fr/hal-01150262}.

\bibitem{Kanenobu:198504}
Taizo Kanenobu.
\newblock {Satellite links with Brunnian properties.}
\newblock {\em Arch. Math.}, 44(4):369--372, 1985.

\bibitem{Kanenobu:1986}
Taizo Kanenobu.
\newblock {Hyperbolic links with Brunnian properties.}
\newblock {\em J. Math. Soc. Japan}, 38:295--308, 1986.

\bibitem{Rolfsen:1976}
Dale Rolfsen.
\newblock {\em Knots and links}.
\newblock {P}ublish or {P}erish, Inc., Houston, 1976, sec. ed. 1990.

\end{thebibliography}

\tableofcontents

\end{document}